# Exact Solution Procedure for the Log-Linear Continuous Knapsack Problem

By

Somdeb Lahiri

(Former Professor) PD Energy University, Gandhinagar (EU-G), India.

ORCID: https://orcid.org/0000-0002-5247-3497

(somdeb.lahiri@gmail.com)

## Abstract

We provide an exact algorithm to solve the log-linear continuous (fractional) knapsack problem. The algorithm is based on two lemmas that follow from the application of weak duality theorem and complementary slackness theorem to the linear optimization problem with linear objective function that is associated with any solution of a linear optimization problem with (differentiable) concave objective function.

## 1. Introduction:

In this paper we provide an exact algorithm to solve the log-linear continuous (fractional) knapsack problem. For continuous knapsack problems, the "main" linear constraint, that is referred to as the "packaging constraint" is appended with the requirement that the decision variables are all less than or equal to one.

There are many approximate algorithms for the general case where the objective functions are non-linear, separable, and concave. However, the only exact algorithms for continuous knapsack problems for non-linear and continuous objective functions are the ones for a "very specific" type of "quasi-linear" objective functions, discussed in Sharkey, Romeijn and Geunes (2011). Log-linear objective functions do not belong to the category of functions discussed by them.

We now provide a motivation for a general continuous knapsack problem, extending the model of a linear continuous knapsack problem, that is available in Dantzig (1957). "Linear extensions" of the continuous linear knapsack problem as formulated by Dantzig, are available in the papers of Sinha and Zoltners (1979) and Witzgal (1980), although the former use it simply for the purpose of applying the "branch and bound algorithm" to the linear programming relaxation of the related integer valued knapsack problem. Such extensions are referred to as "linear continuous multi-choice knapsack problems". However, the undoubted ancestors of continuous knapsack problems of any variety are "Value and Capital" (1939) or for that matter "A Revision of Demand Theory" (1956) both by Sir John Hicks.

For some positive integer 'L', let $\{1, …, L\}$ be a non-empty finite set of activities.

Let m > 0 be the total quantity of a single resource available for allocation among the 'L' activities.

For j∈{1, ..., L}, let $p_j$ > 0 be the maximum amount of the resource that can be allocated to activity j.

Let p = ($p_1$, ..., $p_L$) denote the n-tuple of upper bounds.

Let R: $\mathbb{R}_+^L \to \mathbb{R}$ be a function such for each y = ($y_1$, ..., $y_n$)∈ $\mathbb{R}_+^L$ with $y_j$ being the amount of the resource allocated to activity 'j' for each j∈{1, ..., n}, R(y) is the reward from the resource allocation vector y.

The problem faced by the decision maker is the following.

Maximize R(y)

Subject to $\sum_{j=1}^{L} y_j \leq m$, $y_j \in [0, p_j]$ for j∈{1, ..., L}.

The version of this problem with separable objective function and without any bounds on the allocation of the resource to the activities has been discussed in Lahiri (2002).

Let $R^p: [0,1]^L \to \mathbb{R}$ be the function such that for all x = ($x_1$, ..., $x_n$)∈ $[0, 1]^L$, $R^p(x) = R(p_1 x_1, ..., p_L x_L)$.

Then the above problem is equivalent to the following **continuous knapsack problem**.

Maximize $R^p(x)$

Subject to $\sum_{j=1}^{L} p_j x_j \leq m$, $x_j \in [0, 1]$ for j∈{1, ..., L}.

The continuous knapsack problem has considerable resemblance with one considered by Barucci and Gazzola (2014) and several references therein.

Given any non-empty subset X of $\mathbb{R}^L$ it is easy to see that the set $\{x \in X | \sum_{j=1}^{L} p_j x_j \leq m, x_j \in [0, 1]$ for j∈{1, ..., L}$\}$ = $\{x \in X | \sum_{j=1}^{L} \alpha p_j x_j \leq \alpha m, x_j \in [0, 1]$ for j∈{1, ..., L}$\}$, for all α > 0.

However, it is quite possible that $R^p(x) = R(p_1 x_1, ..., p_L x_L) \neq R(\alpha p_L x_1, ..., \alpha p_L x_L) = R^{\alpha p}(x)$, for some α > 0, and hence the set of solutions to the above maximization problem may be different from the set of solutions to the following maximization problem for some α > 0, α ≠ 1.

Maximize R($\alpha p_1 x_1$, ..., $\alpha p_L x_L$)

Subject to $\sum_{j=1}^{L} \alpha p_j x_j \leq \alpha m$, $x_j \in [0, 1]$ for j∈{1, ..., L}.

In the next section, we introduce the log-linear knapsack problem.

## 2. The log-linear knapsack problem:

Let $\Delta^{L-1} = \{x \in \mathbb{R}_+^L | \sum_{j=1}^{L} x_j = 1\}$ and let $\alpha \in \Delta^{L-1} \cap \mathbb{R}_{++}^L$, $p \in \mathbb{R}_{++}^L$ and m > 0.

The **log-linear knapsack problem** (LLKP) is the following optimization problem:

Maximize $\sum_{j=1}^{L} \alpha_j \ln x_j$

Subject to $\sum_{j=1}^{L} p_j x_j \leq m$, $x_j \leq 1$, $x_j \geq 0$ for all $j \in \{1, \ldots, L\}$.

As noted in the previous section the objective function depends on 'p', but since 'p' is fixed its contribution to the objective function in this case is the constant $\sum_{j=1}^{L} \alpha_j \ln p_j$, which may be ignored in the "ordinally invariant" analysis that follows.

It is easy to see that if x solves the problem, then $x_j > 0$ for all $j \in \{1, \ldots, L\}$. Further, by the strict concavity of the objective function, the optimal solution must be unique.

It is easy to see that if $\sum_{j=1}^{L} p_j \leq m$, then clearly the unique optimal solution for LLKP has all coordinates equal to 1. Hence, let us consider the non-trivial case where $\sum_{j=1}^{L} p_j > m$.

Since, $\sum_{j=1}^{L} p_j > m$, at the optimal solution, the "**packaging constraint**" $\sum_{j=1}^{L} p_j x_j \leq m$, must be satisfied as a strict equality.

Since the packaging constraint must be satisfied with equality, LLKP reduces to the following optimization problem:

Maximize $\sum_{j=1}^{L} \alpha_j \ln x_j$

Subject to $\sum_{j=1}^{L} p_j x_j = m$, $x_j \leq 1$, $x_j \geq 0$ for all $j \in \{1, \ldots, L\}$.

If $\frac{\alpha_j m}{p_j} \leq 1$ for all $j \in \{1, \ldots, L\}$, $x_j^* = \frac{\alpha_j m}{p_j}$ for all $j \in \{1, \ldots, L\}$ is the optimal solution to LLKP, since it is well known that it solves the following maximization problem, for which the variables are <u>not</u> bounded from above by 1, i.e.,

Maximize $\sum_{j=1}^{L} \alpha_j \ln x_j$

Subject to $\sum_{j=1}^{L} p_j x_j \leq m$, $x_j \geq 0$ for all $j \in \{1, \ldots, L\}$.

Hence, suppose that $\{j | \frac{\alpha_j m}{p_j} > 1\} \neq \phi$.

It is well known (see theorem 3.1 in Lahiri (2006) or the theorem in Lahiri (2024)), that $x^* \in \mathbb{R}_{++}^L$ solves LLKP if and only if $x^*$ solves the following linear programming problem denoted Lin-LLKP.

Maximize $\sum_{j=1}^{L} \frac{\alpha_j}{x_j^*} x_j$

Subject to $\sum_{j=1}^{L} p_j x_j = m$, $x_j \leq 1$, $x_j \geq 0$ for all $j \in \{1, \ldots, L\}$.

The dual of this problem is

Minimize $m\lambda + \sum_{j=1}^{L} \mu_j$

Subject to $\lambda p_j + \mu_j \geq \frac{\alpha_j}{x_j^*}$ for all $j \in \{1, \ldots, L\}$, $\mu_j \geq 0$ for all $j \in \{1, \ldots, L\}$.

$\lambda$ is unconstrained in sign.

**3. The associated linear programming problem and its dual:**

If $x^*$ solves Lin-LLKP, then by the Strong Duality theorem of LP, the optimal value of the dual is equal to the optimal value of the primal which in turn is equal to 1.

Suppose $\lambda^*, \mu^*$ solves the dual.

The Complementary Slackness Conditions are the following:

$\mu_j^*(1- x_j^*) = 0$ and $(\lambda^* p_j + \mu_j^* - \frac{\alpha_j}{x_j^*}) x_j^* = 0$ for all $j \in \{1, \ldots, L\}$.

Since we require $x_j^* > 0$, for all $j \in \{1, \ldots, L\}$ we therefore require that $\lambda^* p_j + \mu_j^* - \frac{\alpha_j}{x_j^*} = 0$ for all $j \in \{1, \ldots, L\}$.

$\mu_j^* > 0$ implies $x_j^* = 1$ and $\frac{\alpha_j}{x_j^*} > \lambda^* p_j$, the latter being equivalent to $\frac{\alpha_j}{p_j x_j^*} > \lambda^*$.

If $\mu_j^* = 0$, then $\frac{\alpha_j}{x_j^*} = \lambda^* p_j$, the latter being equivalent to $\frac{\alpha_j}{p_j x_j^*} = \lambda^*$

If $\mu_j^* > 0$ for all $j \in \{1, \ldots, L\}$ then we require $x_j^* = 1$ for all $j \in \{1, \ldots, L\}$.

Thus, $\sum_{j=1}^{L} p_j x_j^* = \sum_{j=1}^{L} p_j > m$, violating the packaging constraint.

Hence $\{j | \mu_j^* = 0\} \neq \phi$. Thus, $\{j | \mu_j^* = 0\} = \{j | \frac{\alpha_j}{p_j x_j^*} = \lambda^* \} \neq \phi$.

**Note:** Since $\alpha_j > 0$ for all $j \in \{1, \ldots, L\}$ and since both $p_j$ and $x_j^* > 0$ by hypothesis for all $j \in \{1, \ldots, L\}$ it must be the case that $\lambda^* > 0$.

Hence, $\lambda^* = \frac{1 - \sum_{j=1}^{L} \mu_j^*}{m} > 0$ would require, $\sum_{j=1}^{L} \mu_j^* < 1$.

If $\mu_j^* = 0$ for all $j \in \{1, \ldots, L\}$, then it must be the case that $\frac{\alpha_j}{x_j^*} = \lambda^* p_j$ for all $j \in \{1, \ldots, L\}$, whence $\alpha_j = \lambda^* p_j x_j^*$ for all $j \in \{1, \ldots, L\}$.

Summing over all $j \in \{1, \ldots, L\}$ and applying the packaging constraint, we get $\lambda^* = \frac{1}{m} \sum_{j=1}^{L} \alpha_j = \frac{1}{m}$, since $\sum_{j=1}^{L} \alpha_j = 1$.

Thus, $x_j^* = \frac{\alpha_j m}{p_j}$ for all $j \in \{1, \ldots, L\}$.

This would require $\frac{\alpha_j m}{p_j} \leq 1$ for all $j \in \{1, \ldots, L\}$, contrary to our assumption that $\{j | \frac{\alpha_j m}{p_j} > 1\} \neq \phi$.

Thus, $\{j | \mu_j^* > 0\} \neq \phi$, whence, $\{j | \frac{\alpha_j}{x_j^*} > \lambda^* p_j \} = \{j | x_j^* = 1\} \neq \phi$.

### 3. Three lemmas:

We now present three lemmas.

**Lemma 1:** If for some $j, k \in \{1, \ldots, L\}$ it is the case that $\frac{\alpha_j}{p_j} \geq \frac{\alpha_k}{p_k}$ and $x_k^* = 1$, the it must be the case that $x_j^* = 1$.

**Proof:** Towards a contradiction suppose $x_j^* < 1$.

Since $x_j^* > 0$, it must be the case that $\lambda^* = \frac{a_j}{p_j x_j^*} > \frac{a_j}{p_j} \geq \frac{a_k}{p_k} = \frac{a_k}{p_k x_k^*} \geq \lambda^*$ implying $\lambda^* > \lambda^*$, which is not possible.

Thus, $x_j^* = 1$. Q.E.D.

An immediate consequence of lemma 1 is the following lemma.

**Lemma 2:** If for some $j \in \underset{k \in \{1,\ldots,L\}}{\operatorname{argmax}} \frac{a_k}{p_k}$ it is the case that $x_j^* = 1$, then <u>for all</u> $j \in \underset{k \in \{1,\ldots,L\}}{\operatorname{argmax}} \frac{a_k}{p_k}$ it must be the case that $x_j^* = 1$.

As a result of lemmas 1 and 2, we have the following lemma.

**Lemma 3:** If $\{j | \frac{a_j m}{p_j} > 1\} \neq \phi$, then $x_j^* = 1$ for all $j \in \underset{k \in \{1,\ldots,L\}}{\operatorname{argmax}} \frac{a_k}{p_k}$. Further, $\underset{k \in \{1,\ldots,L\}}{\operatorname{argmax}} \frac{a_k}{p_k}$ is a non-empty proper subset of $\{1, \ldots, L\}$.

**Proof:** Suppose $\{j | \frac{a_j m}{p_j} > 1\} \neq \phi$ and towards a contradiction suppose that for some $j \in \underset{k \in \{1,\ldots,L\}}{\operatorname{argmax}} \frac{a_k}{p_k}$ it is the case that $x_j^* < 1$. Then, by lemma 2, it must be the case that $x_j^* < 1$, <u>for all</u> $j \in \underset{k \in \{1,\ldots,L\}}{\operatorname{argmax}} \frac{a_k}{p_k}$.

By lemma 1, it follows that $x_j^* < 1$ for all $j \in \{1, \ldots, L\}$ and hence $\mu_j^* = 0$ for all $j \in \{1, \ldots, L\}$.

Thus, $\{j | \mu_j^* > 0\} = \phi$, leading to a contradiction.

Thus, it must be the case that $x_j^* = 1$ for all $j \in \underset{k \in \{1,\ldots,L\}}{\operatorname{argmax}} \frac{a_k}{p_k}$.

Since we have assumed, $\sum_{j=1}^{L} p_j > m$, $x_j^* = 1$ for all $j \in \{1, \ldots, L\}$, would lead to a violation of the packaging constraint.

Hence, it must be the case that $\underset{k \in \{1,\ldots,L\}}{\operatorname{argmax}} \frac{a_k}{p_k}$ is a proper subset of $\{1, \ldots, L\}$.

This, proves the lemma. Q.E.D.

### 4. The "exact" algorithm:

Lemma 3 yields an iterative method of solving LLKP.

**Step 1:** If $\frac{a_j m}{p_j} \leq 1$ for all $j \in \{1, \ldots, L\}$, $x_j^* = \frac{a_j m}{p_j}$ for all $j \in \{1, \ldots, L\}$ is the optimal solution to LLKP. If not, i.e., $\{j | \frac{a_j m}{p_j} > 1\} \neq \phi$, then let $x_j^* = 1$ for all $j \in \underset{k \in \{1,\ldots,L\}}{\operatorname{argmax}} \frac{a_k}{p_k}$.

Let $J^1 = \{1, \ldots, L\} \setminus \underset{k \in \{1,\ldots,L\}}{\operatorname{argmax}} \frac{a_k}{p_k}$.

Since $\sum_{j=1}^{L} p_j > m$, it must be the case that $\sum_{j \in J^1} p_j > m - \sum_{j \in \underset{k \in \{1,\ldots,L\}}{\operatorname{argmax}} \frac{a_k}{p_k}} p_j = m - \sum_{j \notin J^1} p_j$

**Step 2:** Consider the revised problem

Maximize $\sum_{j \in J^1} \alpha_j \ln x_j$

Subject to $\sum_{j \in J^1} p_j x_j = m - \sum_{j \notin J^1} p_j$,  $x_j \leq 1$, $x_j \geq 0$ for all $j \in J^1$.

The above problem is equivalent to the following problem.

Maximize $\sum_{j \in J^1} \alpha_j^{(1)} \ln x_j$

Subject to $\sum_{j \in J^1} p_j x_j = m - \sum_{j \notin J^1} p_j$,  $x_j \leq 1$, $x_j \geq 0$ for all $j \in J^1$, where for all $h \in J^1$, $\alpha_h^{(1)} = \dfrac{\alpha_h}{\sum_{j \in J^1} \alpha_j}$.

If $\dfrac{\alpha_j^{(1)}(m - \sum_{j \notin J^1} p_j)}{p_j} \leq 1$ for all $j \in J^1$, then let $x_j^* = \dfrac{\alpha_j^{(1)}(m - \sum_{j \notin J^1} p_j)}{p_j}$ for all $j \in J^1$. If not, i.e. $\{j \in J^1 \mid \dfrac{\alpha_j^{(1)}(m - \sum_{j \notin J^1} p_j)}{p_j} > 1\} \neq \emptyset$, then let $x_j^* = 1$ for all $j \in \underset{k \in J^1}{\operatorname{argmax}} \dfrac{\alpha_k}{p_k}$.

Since this process cannot go on forever, we will finally arrive at the situation where either all the values of $x_j^*$ for $j \in \{1, \ldots, L\}$ have been obtained by repeated application of the above procedure or there exists $h \in \{1, \ldots, L\}$ such that the values $x_j^*$ for all $j \in \{1, \ldots, L\} \setminus \{h\}$ have been determined by repeated application of the above procedure, so that $x_j^* = 1$ for all $j \in \{1, \ldots, L\} \setminus \{h\}$. In the latter case, $x_h^* = \dfrac{m - \sum_{j \neq h} p_j}{p_h}$.

Since by hypothesis, $\sum_{j=1}^{l} p_j > m$ the only way that by repeated application of the above procedure $x_j^*$ would have been determined for all $j \in \{1, \ldots, L\}$ without resorting to the formula $x_h^* = \dfrac{m - \sum_{j \neq h} p_j}{p_h}$ for some $h \in \{1, \ldots, L\}$, is when the solution to the last LLKP problem yields the solution for more than one variable and agrees with the one obtained by dropping the upper bound of 1 on the variables whose values are determined in the last LLKP, with all values of the variables determined in previous LLKP's being equal to 1.

## 5. A numerical example to illustrate the procedure:

Let $L = 3$ and suppose $\alpha_j = \dfrac{1}{3}$ for $j = 1, 2, 3$. Let $p_1 = 1$, $p_2 = 2$, $p_3 = 3$ and $m = 5$.

Thus, $\dfrac{\alpha_1}{p_1} = \dfrac{1}{3}$, $\dfrac{\alpha_2}{p_2} = \dfrac{1}{6}$ and $\dfrac{\alpha_3}{p_3} = \dfrac{1}{9}$.

The log-linear knapsack problem is the following.

Maximize $\dfrac{1}{3} \ln x_1 + \dfrac{1}{3} \ln x_2 + \dfrac{1}{3} \ln x_3$,

Subject to $x_1 + 2x_2 + 3x_3 \leq 5$, $1 \geq x_j \geq 0$ for $j = 1, 2, 3$.

Note that $p_1 + p_2 + p_3 = 6 > 5 = m$. Hence $(1, 1, 1)$ is not feasible.

For the unconstrained problem,

Maximize $\dfrac{1}{3} \ln x_1 + \dfrac{1}{3} \ln x_2 + \dfrac{1}{3} \ln x_3$,

Subject to $x_1 + 2x_2 + 3x_3 \leq 5$, $x_j \geq 0$ for $j = 1, 2, 3$,

the optimal solution is $(\frac{5}{3}, \frac{5}{6}, \frac{5}{9})$ which is not feasible for this log-linear knapsack problem.

Since, $\frac{\alpha_1}{p_1} > \frac{\alpha_j}{p_j}$ for $j = 2, 3$, according to our algorithm, at the optimal solution $x_1^* = 1$.

Since $\frac{\frac{1}{3}}{\frac{1}{3} + \frac{1}{3}} = \frac{1}{2}$, we consider the revised problem

Maximize $\frac{1}{2}\ln x_2 + \frac{1}{2}\ln x_3$,

Subject to $2x_2 + 3x_3 \leq 5$, $1 \geq x_j \geq 0$ for $j = 2, 3$.

For the unconstrained problem,

Maximize $\frac{1}{2}\ln x_2 + \frac{1}{2}\ln x_3$,

Subject to $2x_2 + 3x_3 \leq 4$, $x_j \geq 0$ for $j = 2, 3$,

the optimal values of $x_2$ and $x_3$ are 1 and $\frac{2}{3}$ respectively and these values are feasible for the log-linear knapsack problem for $x_2$ and $x_3$.

Hence, these must be the optimal values of $x_2$ and $x_3$ for the log-linear knapsack problem we started out with. Thus, the optimal solution for the original log-linear knapsack problem is $x_1^* = 1$, $x_2^* = 1$ and $x_3^* = \frac{2}{3}$.

**Acknowledgment:** My interest in the general version of continuous knapsack problems formulated in the introductory section of the paper, was "triggered" by a search for a meaningful interpretation of "prices in the utility function" that is assumed by Amit Goyal in a short note that he had prepared on "Veblen" goods (https://economics.stackexchange.com/questions/49928/upward-sloping-demand-curves-can-t-exist/54899#54899). I would like to thank him for both the note as well as the discussion that we had based on it, though I am unable to find any convincing "economic" reason for the violation of "law of demand" that his analysis implies. I would also like to thank Seyyed Ahmad Adalat Panah for his observations and for suggesting the inclusion of a numerical example, which resulted in section 5 of the paper. Further, I would like to thank Samir Kumar Neogy for his encouraging response about the contents of this paper. However, no one other than the author is responsible for residual errors.